\documentclass[11 pt]{article}
\usepackage{makeidx,graphics}
\usepackage{graphicx}
\usepackage{color}
\usepackage{hyperref}
\usepackage{amsmath,amsthm,amssymb, amscd}
\usepackage{geometry}
\color[cmyk]{0,0,0,1}

\DeclareSymbolFont{AMSb}{U}{msb}{m}{n}
\DeclareMathSymbol{\N}{\mathbin}{AMSb}{"4E}
\DeclareMathSymbol{\Z}{\mathbin}{AMSb}{"5A}
\DeclareMathSymbol{\R}{\mathbin}{AMSb}{"52}
\DeclareMathSymbol{\Q}{\mathbin}{AMSb}{"51}
\DeclareMathSymbol{\I}{\mathbin}{AMSb}{"49}
\DeclareMathSymbol{\C}{\mathbin}{AMSb}{"43}

\newcommand{\Ad}{{\rm Ad}}

\newcommand{\F}{\mathbb F}
\newcommand{\quat}{\mathbb H}

\newcommand{\tr}{{\rm tr}}

\title{On 3 and 4 dimensional regular solids\\ Part 1: The 4-simplex generates the free group}

\author{Adrian Ocneanu}

\begin{document}
\date{}
\maketitle

\begin{abstract}
The 4-simplex has vertices 5 unit quaternions, which we arrange so that one of them is the unit. We show that the remaining vertices are the generators of a free group. 

For the proof, we introduce a new alternating length on words in free groups. We show that for words in simplex vertices the necklace form of the alternating length can be read number theoretically, as the logarithm of the algebraic denominator of their trace.

\end{abstract}

\section{Main results}

\newtheorem{thm}{Theorem}
\newtheorem{prop}{Proposition}
\newtheorem{rem}{Remark}

We choose the vertices of the regular 4D simplex as unit quaternions, with one vertex the unit 1.
\begin{eqnarray*}
q_0&=&(1,0,0,0)\\
q_1&=&\left ( -\frac{1} {4},\frac{\sqrt 5} {4},\frac{\sqrt 5} {4},\frac{\sqrt 5} {4} \right )\\
q_2&=&\left ( -\frac{1} {4},-\frac{\sqrt 5} {4},-\frac{\sqrt 5} {4},\frac{\sqrt 5} {4} \right )\\
q_3&=&\left ( -\frac{1} {4},\frac{\sqrt 5} {4},-\frac{\sqrt 5} {4},-\frac{\sqrt 5} {4} \right )\\
q_4&=&\left ( -\frac{1} {4},-\frac{\sqrt 5} {4},\frac{\sqrt 5} {4},-\frac{\sqrt 5} {4} \right )
\end{eqnarray*}

The order of the generators is chosen so that the volume of the simplex is positive, a fact which we shall use later.

Our main result is the following.
\begin{thm}[Main Theorem]
The vertices $q_1,q_2,q_3,q_4$ are generators of a free group.
\end{thm}

\section{Alternating length and clutched necklaces in free groups}

For a word $w$ in a free group $F$ we define the \textbf{alternating norm} $|w|_{\rm alt}$ as follows. The word $w$ has a unique  \textbf{alternating pairs expression}
$$w=(a_1 a_2^{-1})(a_3 a_4^{-1})\dots (a_{2n-1} a_{2n}^{-1}) $$
with consecutive terms $ 1\ne a_1\ne a_2\ne\dots\ne a_{2n-1} \ne a_{2n}$ and either all $a_i\in\{1,g_1,\dots,g_m\}$ or all $a_i\in\{1,g_1^{-1},\dots,g_m^{-1}\}$. Define $|w|_{\rm alt}=n$, the number of such pairs. Note that $|w|_{\rm alt}=0$ iff $w=1$.

\begin{prop} 
For a word $w$ in reduced form, with length $|w|$ and having $s$ sign changes between the powers of its generators, we have  $|w|_{\rm alt}=|w|- \lceil s/2\rceil $.
\end{prop}

For instance $w=g_1^{-1}g_2 g_3 g_1^{-1} = (g_1^{-1}g_2)(1^{-1} g_3)(g_1^{-1} 1)\in \F_3$ has $s=2$ and $|w|_{\rm alt}=3 = |w|-\lceil s/2\rceil = 4-1$.

Due to the symmetry of the definitions we may assume without loss of generality that $w$ begins with a generator $g^{+ 1}$. We prove the statement by induction on $|w|$.

For the inductive proof, assume first that the number $s$ of sign changes of $w$ is odd. In this case $\lceil (s+1)/2\rceil = \lceil s/2\rceil.$
Since $s$ is odd, $w$ will end with a generator $h^{- 1}$, and have the last alternating expression pair $(\bullet^{+1}h^{- 1})$. By induction the number of alternating parentheses of $w$ is $|w|- \lceil s/2\rceil$. An extra generator $t^{\pm 1}$ increases the word length by 1, increases the number of alternating pairs by 1, while $\lceil s/2\rceil $ stays the same, so the equality is preserved. Assume now that $s$ is even. Then the alternating form of $w$ ends with $(h^{+1} 1^{-1})$. An extra generator  $t^{\pm 1}$ increases the word length by 1. An additional generator $t^{-1}$ at the end makes the last parenthesis $(h^{+1} t^{-1})$ and increases the number of sign changes, and $\lceil s/2\rceil$, by 1. An extra generator $t^{+1}$ increases the number of parentheses by 1 but keeps the number of sign changes the same. In both cases, the equality is preserved.

Assume that the reduced word $w$ is shortest among its conjugates. That means that if $w$ starts with a generator $g^{\pm 1}$ then $w$ does not end with $g^{\mp 1}$. Equivalently, $w$ is {\bf necklace reduced}, that is, if we put $w$ in {\bf necklace form} on a circle then there are no cancelations.

In the necklace form of a necklace reduced word $w$, the total number of sign changes $s'$ is even, $s'=2\lceil s/2\rceil $, and in that case $|w|_{\rm alt}=|w|- s'/2$. Remark that in this case moving the first letter of $w=g^{\pm 1}v$ from front to back to $w'=vg^{\pm 1}$  gives a necklace reduced word with the same alternating length.

We shall need a slightly more general construction, which we call the {\bf clutch necklace}. Suppose that we are given a permutation $\pi$ of the generators of the free group. Extend $\pi$ to words. A reduced word $w$ is called $\pi$-reduced if it is shortest among all its $\pi$-conjugates $u w\pi(u)^{-1}$. A reduced word $w$ which starts with the generator $g^{\pm 1}$ is $\pi$-reduced if it does not end with  $\pi(g)^{\mp 1}$. In this case $w$ can be considered to be on a necklace which has at the connecting point between the end and the beginning of $w$ a clutch marked with $\pi$. A word $u$ jumping over the clutch, from the beginning to the end of $w$, is replaced by $\pi(u)$. Moving a beginning portion $u$ of a $\pi$-reduced $w=uv$ with $|w|= |u|+|v|$ over the clutch to the end to get $w' = v\pi(u)$ keeps $w'$ $\pi$-reduced and with $|w'|_{\rm alt} = |w|_{\rm alt}$.

The  {\bf $\pi$-reduced} form of $w$, ${\rm red}_\pi(w)$, is obtained from the reduced word $w$ by canceling a generator $g^{\pm 1}$ at the beginning of $w$ with a generator $\pi(g)^{\mp 1}$, repeated as long as this is possible.

The length which we shall need is the \textbf{reduced $\pi$-alternating norm} $||w||_\pi = |{\rm red}_\pi(w)|_{\rm alt}$. 

Note that the not necessarily trivial word $w= g\pi(g)^{-1}$ has ${\rm red}_\pi(w) = 1$ and so $||w||_\pi = 0$.

\begin{rem} 
We shall introduce in a second paper in this series the {\bf twin dodecahedra} with vertices stereographic projections of quotients $q_iq_j^{-1}$ and $q_i^{-1}q_j$. The alternating length $ |\cdot |_{\rm alt}$ as well as the lengths $||\cdot ||_\pi$ are geometrically related to these vertices.
\end{rem}

\section{Quaternionic permutations}

Let us return now to the case of the regular 4-simplex vertices $q_0=1,q_1,q_2,q_3,q_4\in\quat$. Let $\F_4$ denote the free group on generators $g_1,g_2,g_3,g_4$ and denote by $\rho:g_i\mapsto q_i\in\quat$ the representation which maps its generators onto the nontrivial simplex vertices. 

We denote by $i_1, i_2, i_3$ the generators of the quaternion group, the $i,j,k$ in the notation of Hamilton.

We have $\Ad_{ i_1}: i_1\mapsto i_1,i_2\mapsto -i_2,i_3\mapsto -i_3$ and similar properties for $i_2$ and $i_3$. As the conjugation by a quaternion generator changes the signs of the other two generators it permutes the vertices $q_1,q_2,q_3,q_4$ in pairs, as in the table below. We denote by $\pi_k$ the permutation of the indices $1,2,3,4$ with $\Ad_{i_k}(q_i) = q_{\pi_k(i)}$. 
\begin{center}
\begin{tabular}{c|c|cccc|c|}
$k$&$i_k$ & $q_1$& $q_2$& $q_3$& $q_4$& $\pi_k\ {\rm cycles}$\\
\hline 
  1&$i_1$& $+$& $-$& $+$& $-$&(1,3)(2,4)\\ 
 2&$i_2$& $+$& $-$& $-$&$+$&(1,4)(2,3)\\ 
 3& $i_3$& $+$& $+$& $-$&$-$&(1,2)(3,4)\\ 
\end{tabular}
\end{center}
The map $k\mapsto \pi_k$ is thus the map of the quaternion group into the Klein group, in its representation by permutations.

We extend $\pi_k$ to products of simplex vertices. For completeness we let $i_0=1$ and $\pi_0={\rm id}$. We also extend the permutation $\pi_k$ to the generators and words in $\F_4$.

Let $k\in\{0,1,2,3\}$. To simplify notation, we shall omit $\pi$ and define for $w\in \F_4$ as we did in the general free group case the \textbf{$k$-reduction} ${\rm red}_{k}(w) = {\rm red}_{\pi_k}(w)$ and the \textbf{reduced $k$-alternating norm}  by $||w||_k=||w||_{\pi_k} = |{\rm red}_{k}(w)|_{\rm alt}|$.

\section{The algebraic denominator and its shifted logarithm {\rm lad}}
Recall that a number $x$ is an \textbf{algebraic integer} if it satisfies an integer equation with the coefficient of the maximal power 1. We call an algebraic integer $x$ \textbf{even} if $x/2$ is an algebraic integer, and otherwise we call it \textbf{odd}. The sum of even algebraic integers is even, so the sum of an odd and an even algebraic integer is odd. The sum of odd algebraic integers may be odd. For instance the  two roots $(1\pm\sqrt 5)/2$ of the golden ratio equation $x^2-x-1=0$, which are odd algebraic integers, have sum 1. 

Remark that words in $q_i^{\pm 1}$ have denominator, and thus algebraic denominator, a power of 2. For an algebraic number $x=y/2^{k+1}$ where $y$ is an odd algebraic integer, denote by ${\rm lad}(x)=k$ the {\bf shifted logarithm of the algebraic denominator}, with 
${\rm lad}(0) = - \infty$. 

We have ${\rm lad}(xy)={\rm lad}(x)+{\rm lad}(y)+1$ (because of our shift) and ${\rm lad}(x+y)\le \max({\rm lad}(x),{\rm lad}(y))$.
From the fact that the sum between an odd and an even algebraic integer is odd, we see that if ${\rm lad}(x) \ne {\rm lad}(y)$ then we have an equality, ${\rm lad}(x+y) = \max({\rm lad}(x),{\rm lad}(y))$.

\section{The Algebraic Denominator Theorem}

Our main theorem has the following precise form.

\begin{thm}[The Algebraic Denominator Theorem]
 Let $w \in \F_4$ and its quaternion representation $x = (x_0,x_1,x_2,x_3) =\rho(w)$.
 
Then $||w||_k=0$ iff ${\rm red}_k(w) = 1$ iff $k=0$ and $x = 1$ or $k>0$ and $x_k =0$.
Otherwise
$${\rm lad}(x_k) = ||w||_k \ge 1$$
\end{thm}

For instance $w=g_1g_2^{-1}=g_1\pi_3(g_1)^{-1}$ is 3-conjugate to $g_1^{-1} g_1 =1 ={\rm red}_3(w)$, so $||w||_3=0$. It has, for $k=0,1,2,3$, $k$-reduced forms ${\rm red}_k(w)$ equal to $w,w,w,1$ respectively. Thus $||w||_k=|{\rm red}_k(w)|_{\rm alt}$ equals $1,1,1,0$ respectively. The representation of $w$ is
$$x=\rho(w)=q_1 q_2^{-1}=\left(\frac{-1}{4},\frac{-5-\sqrt 5}{8},\frac{5-\sqrt 5}{8},0\right).$$
Its components $x_k$ have algebraic denominators $4,4,4,0$ which have $-1+\log_2($alg.denom.($x_k$)) = lad($x_k$) equal to 1,1,1,$-\infty$ respectively.

Similarly $w=g_1g_2$ is $k$-reduced for all $k$ and as $w=(g_11^{-1})(g_21^{-1})$ has $|w|_{\rm alt}=2$ and so $||w||_k=2$ for all $k$. As
$$x=\rho(w)=q_1 q_2=\left(\frac{3}{8},\frac{3}{8},\frac{3}{8},\frac{-5}{8}\right)$$ 
we have $-1+\log_2($alg.denom.($x_k$)) = lad($x_k$)=2=$||w||_k$ for all $k$. The identity ${\rm lad}(x_k) = ||w||_k \ge 1$ if $||w||_k \ge 1$ is thus immediate for words $w$ of length 1 or 2.

\section{The Proof}

We now start the proof of the theorem. Let $w\in F_4$. We have remarked before that the only word with alternating length 0 is 1, so $||w||_k=|{\rm red}_k(w)|_{\rm alt} = 0$ iff  ${\rm red}_k(w) = 1$. Note that ${\rm red}_0(w)$ is conjugate to $w$, so if ${\rm red}_0(w) = 1$ then $w$ =1, and $x=\rho(w)$ =1.

Let $k\in\{1,2,3\}$ and let $u\in \F_4$ be the word with 
${\rm red}_k(w) = u w \pi_k(u)^{-1}$. Then 
\begin{eqnarray*}
x_k&=&-\tr(i_k x)=-\tr(i_k \rho(w)) = -\tr(i_k \rho(u)^{-1}\rho(u)  \rho(w))= -\tr(\rho(\pi_k(u))^{-1}i_k \rho(u)  \rho(w)) =\\
&=& -\tr(i_k \rho(u w \pi_k(u)^{-1})) = -\tr(i_k \rho({\rm red}_k(w))).
\end{eqnarray*}
Thus if  ${\rm red}_k(w) = 1$ then $x_k=\tr(i_k)=0$.

This proves one direction of the first assertion of the theorem.

The opposite direction implies that if $x=\rho(w) =1$ then ${\rm red}_0(w)=1$, so it shows that $\rho$ is injective and its image is free. 

It will be implied by the main equality, which states that if $w\in \F_4$ is $k$-reduced and has length $|w|\ge 1$, so that $||w||_k =|w|_{\rm alt} \ge 1$,  then ${\rm lad}(x_k) = ||w||_k \ge 1$. In particular if $k=0$ and $|w|>1$ then ${\rm lad}(x_0) \ge 1$ so $x\ne 1$.

We now proceed to prove the main equality by induction on the length $|w|$ of the reduced word $w\in \F4$. Note that for $|w|=1$ we have $w=g_i^{\pm 1},$ so $x=\rho(w) = q_i^{\pm 1}$ has all components $-1/4$ or $\pm\sqrt 5/4$ with algebraic denominator 4, and thus for all $k$ we have $||w||_k=1=\log_2(4)-1={\rm lad}(x_k).$

We shall make use of the following elementary quadratic identities, together with their analogs with permuted indices, and with all $q_i$ replaced by $q_i^{-1}$. Remark that the problem is symmetric in the generators $q_1,\dots,q_4$. Moreover the Galois symmetry $\sqrt 5 \mapsto -\sqrt 5$ maps $q_i$ into $q_i^{-1}$.
\begin{eqnarray*}
q_1^2&=&-\frac 1 2 q_1-1\\
q_1 q_2&=&\frac {\sqrt 5} 2 q_1 i_3+1=\frac {\sqrt 5} 2  i_3 q_2 +1\\
q_1 q_2^{-1}q_1&=& -\frac 1 2 q_1-q_2 \\
q_1 q_2^{-1}q_3&=&\frac {5+\sqrt 5} 4 q_1 i_1+q_4=\frac {5+\sqrt 5} 4 i_1 q_3+q_4  
\end{eqnarray*}

All the equalities above have the form $$q_1q_i^{-1}q_j=c\  q_1 i_k+q_l=c\ i_k q_j+q_l$$ where $c$ is a number with algebraic denominator 2, $k>0$ (i.e. $q_k\ne 1$) and we have $j>0$ iff $l>0$.

If $w\in \F_4$ is a 0-reduced word $w={\rm red}_0(w)$ which starts with the left hand side $w= g_1g_i^{-1}g_j v$ as a word in $\F_4$, then $v$ does not end in $g_1^{-1}$. 

Then $u = g_j v$ is $k$-reduced, again since $v$ does not end in $g_1^{-1}$ and  $u=g_j v$ is $k$-conjugate to $v\pi_k(g_j)=v g_1$.
Thus $$||w||_0=|w|_{\rm alt}=|g_1g_i^{-1}g_j v|_{\rm alt}=1+|g_j v|_{\rm alt}=1+||u||_k.$$
The word corresponding to the last term $r=g_l v$ has an alternating form strictly shorter than $w$, so $||w||_0>||r||_0.$

We multiply the equality with $\rho(v)$ to obtain
$$\rho(w) = q_1q_i^{-1}q_j\rho(v)=c\  i_k q_j\rho(v)+q_l\rho(v) =c\ i_k \rho(u) + \rho(r)$$
and we take the trace. Notice that $\tr(i_k \rho(u)) = -(\rho(u))_k$, the $k$-th component of the quaternion $\rho(u)$.
By induction we have the relation between shifted algebraic denominators
$${\rm lad}(\tr(c\ i_k\rho(u)))= 1+{\rm lad}(\tr(i_k\rho(u))) = 1+{\rm lad}((\rho(u))_k) = 1+||u||_k = ||w||_0$$
while 
$${\rm lad}(\tr(\rho(r))) =||r||_0 < ||w||_0$$
Since the two terms above have unequal algebraic denominators, the algebraic denominator of their sum is the maximum of the two, so 
$${\rm lad}(\tr(\rho(w))) ={\rm lad}(\tr(c\ i_k \rho(u))) + {\rm lad}(\tr( \rho(r))) =  {\rm lad}(\tr(c\ i_k \rho(u))) = ||w||_0$$
which proves the equality by induction.

Now a 0-reduced word $w$ of length $\ge 3$ has either, up to conjugacy, two consecutive $g_s g_t$ or $g_s^{-1} g_t^{-1}$ with $s,t>0$, case which is covered by the first two equalities above (or the similar ones with negative powers of the generators), or else it has an alternating powers form covered by the last 2 equalities. The case of the norm $||w||_k$ with $k>0$ is entirely similar and left to the reader.

Up to the manifest symmetry in permuting the generators, inverting them and changing the index $k$, we have shown that  $||w||_k = {\rm lad}(x_k)$ by induction. 

This ends the proof of the structure theorem, and proves that the subgroup of $\quat$ generated by the nontrivial vertices of the 4-simplex is free.


\section{Acknowledgments}
We thank Nick Early who recognized the importance of our result and strongly encouraged us to publish it. We obtained the result in 2006 as we constructed a large scale stainless steel 4-dimensional sculpture, the Octacube,  at Penn State, and we lectured on 3 and 4 dimensional polytopes and on the quantum subgroups of $SU(2)$, which we defined, in an advanced undergraduate and respectively a graduate course at Penn State.
\end{document}